\documentclass[10pt]{amsart}

\usepackage[latin1]{inputenc}
\usepackage{amsfonts}
\usepackage{amsmath}
\usepackage{amssymb}
\usepackage{amsthm}
\usepackage[margin=1in]{geometry}
\usepackage{graphicx}

\theoremstyle{plain}
	\newtheorem{theorem}{Theorem}

	\newtheorem*{theorem*}{Theorem}
	\newtheorem*{lemma*}{Lemma}
	\newtheorem*{claim*}{Claim}
	\newtheorem*{prop*}{Proposition}
	
\theoremstyle{definition}
	
	\newtheorem*{question*}{Question}
	\newtheorem*{note*}{Note}

\newcommand{\Z}{\mathbb{Z}}
\newcommand{\R}{\mathbb{R}}
\newcommand{\Q}{\mathbb{Q}}
\newcommand{\C}{\mathbb{C}}

  \def\Dbar{\leavevmode\lower.6ex\hbox to 0pt{\hskip-.23ex \accent"16\hss}D}
  \def\cftil#1{\ifmmode\setbox7\hbox{$\accent"5E#1$}\else
  \setbox7\hbox{\accent"5E#1}\penalty 10000\relax\fi\raise 1\ht7
  \hbox{\lower1.15ex\hbox to 1\wd7{\hss\accent"7E\hss}}\penalty 10000
  \hskip-1\wd7\penalty 10000\box7}
  \def\cfudot#1{\ifmmode\setbox7\hbox{$\accent"5E#1$}\else
  \setbox7\hbox{\accent"5E#1}\penalty 10000\relax\fi\raise 1\ht7
  \hbox{\raise.1ex\hbox to 1\wd7{\hss.\hss}}\penalty 10000 \hskip-1\wd7\penalty
  10000\box7}

\numberwithin{equation}{section}

\begin{document}

\title{Branching of Automorphic Fundamental Solutions}
\author{Amy T. DeCelles}
\address{Department of Mathematics, University of St. Thomas, 2115 Summit Ave., St. Paul, MN 55105}
\email{adecelles@stthomas.edu}
\urladdr{http://personal.stthomas.edu/dece4515}
% 11F72 Discontinuous groups and automorphic forms: spectral theory; Selberg trace fmla
%11F55  Discontinuous groups and automorphic forms: other groups and their modular and automorphic forms (several variables)
%11F41 Automorphic forms on Hilbert modular surfaces
\thanks{The author thanks Paul Garrett, for posing this problem, Jeffrey Lagarias, for suggesting further work on the problem, and both of them for helpful questions and conversations.  The author was partially supported by a research grant from the University of St. Thomas.}

\begin{abstract}
Automorphic fundamental solutions and, more generally, solutions of automorphic differential equations, play a key role in the Diaconu-Garrett-Goldfeld prescription for spectral identities involving moments of L-functions \cite{diaconu-garrett09, diaconu-garrett10, diaconu-garrett-goldfeld10a} as well as other applications, including an explicit formula relating the number of lattice points in a symmetric space to the automorphic spectrum \cite{decelles-lattice2011}.  In this paper we discuss two cases in which the automorphic fundamental solution exhibits branching: pathwise meromorphic continuations may differ by a term involving an Eisenstein series. 
\end{abstract}

 \maketitle

 \section{Introduction}
 
Solutions of automorphic differential equations underlie the Diaconu-Garrett-Goldfeld prescription for spectral identities involving second moments for arbitrary Rankin-Selberg integral representations of L-functions \cite{diaconu-garrett-goldfeld10a}.  This prescription is a vast generalization of the constructions of moment identities in their earlier papers, from which they extracted subconvex bounds for $GL_2$ automorphic L-functions \cite{diaconu-goldfeld06, diaconu-goldfeld07, diaconu-garrett09, diaconu-garrett10}.  Essential to their prescription is a Poincar\'{e} series, whose data was originally constructed in imitation of Good's kernel \cite{good82}; characterizing the Poincar\'{e} series as the solution to an automorphic differential equation allows generalization from $GL_2$ to higher rank.  The automorphic spectral expansion of such a Poincar\'{e} series is heuristically immediate and can be legitimized using automorphic Sobolev theory, developed in \cite{decelles-lattice2011}.  In general, explicit geometric expressions for solutions of automorphic differential equations are very difficult to obtain; see \cite{decelles-designedPe2013} for some examples, including the automorphic fundamental solution that is used in the lattice-point counting application in \cite{decelles-lattice2011} and is suitable for constructing moment identities for $GL_n(\C) \times GL_n(\C)$ Rankin-Selberg L-functions.  Superficially, the spectral expansion of the automorphic fundamental solution appears to be invariant under a transformation of an auxiliary complex parameter $w$, but a closer look reveals, in certain cases, \emph{branching} in $w$, eliminating the possibility of a straightforward functional equation.

%While the spectral expansion of the automorphic fundamental solution suggests invariance under $w \to 1-w$, where $w$ is a complex ``eigenvalue'' parameter, pathwise meromorphic continuations past the line $\mathrm{Re}(w) = \tfrac{1}{2}$ along different $w$-paths may differ by a term involving an Eisenstein series.  In this paper, we describe two cases in which this occurs: the Hilbert-Maass case and the $GL_3$ case.
 
%%%%%%
 
Let $G$ be a semi-simple Lie group, $K$ its maximal compact subgroup, and $\Gamma$ a discrete subgroup.  Consider the solution of the following differential equation on the arithmetic quotient $X = \Gamma \backslash G/K$:
%$$(\Delta - \lambda_z)^{\nu} \; v^{\mathrm{afc}}_z \;\; = \;\; \delta_{x_o}$$
$$(\Delta - \lambda_w)^{\nu} \; u_w \;\; = \;\; \delta_{z_o},$$
where the Laplacian $\Delta$ is the image of the Casimir operator for $\mathfrak{g}$, $\lambda_w$ a complex parameter, $\nu$ an integer, and $\delta_{z_o} = \delta_{\Gamma \cdot 1\cdot K}$ the Dirac delta distribution at the basepoint in $\Gamma \backslash G/K$.  In rank one, we parametrize $\lambda_w$ as $\lambda_w = w(w-1)$, and in higher rank, $\lambda_w = w^2 - |\rho|^2$, where $\rho$ is the half-sum of positive roots.  We recall the following results from \cite{decelles-lattice2011}.  Global automorphic Sobolev theory ensures the existence of a solution $u_w$, unique in global automorphic Sobolev spaces.  The solution has a transparent automorphic spectral expansion, converging in a global automorphic Sobolev space, for $\mathrm{Re}(w)$ sufficiently large.   Further, a global Sobolev embedding theorem ensures that, by choosing $\nu$ sufficiently large, we may ensure that the spectral expansion converges uniformly pointwise or in any $C^k$-topology that we wish. 

Interestingly, the fundamental solution may exhibit \emph{branching} in the complex variable $w$: meromorphic continuations along different $w$-paths in the complex plane may differ by a term of moderate growth.  In particular, the resulting function may lie outside of global automorphic Sobolev spaces.  We discuss two such cases below.%: when the groundfield is a totally real number field of degree $n>1$ over $\Q$ and the case of the $GL_3(\R)$ automorphic fundamental solution.

Branching of fundamental solutions on symmetric spaces has been discussed by Mazzeo and collaborators in several papers; see, e.g. \cite{mazzeo-vasy}, and by Strohmaier \cite{strohmaier}.  See \cite{bruin}, for a discussion of automorphic Green functions on $\Gamma \backslash \mathfrak{H}$, where $\Gamma$ is a Fuchsian group and \cite{oda-tsuzuki}, for automorphic Green functions  with logarithmic singularities along modular divisors in a modular variety.

\section{Branching of Hilbert-Maass Fundamental Solutions}
Recall the general set-up from the introduction, which invokes \cite{decelles-lattice2011} for the existence and uniqueness of solutions to automorphic differential equations and convergence of their spectral expansions in global automorphic Sobolev spaces.  For complete descriptions of automorphic spectral expansions see \cite{langlands76,moeglin-waldspurger95}.

Let $k$ be a totally real number field of degree $n>1$ over $\Q$ with $\mathfrak{o}$ its ring of integers.  For simplicity, suppose $\mathfrak{o}$ has narrow class number one, so that $SL_2(\mathfrak{o})$ is unicuspidal.  Let  $\sigma_1, \dots \sigma_n$ be the archimedean places of $k$, and let $SL_2(\mathfrak{o})$ act on $\mathfrak{H}^n$, componentwise, as usual: for $\gamma \in SL_2(\mathfrak{o})$ and $z = (z_1, \dots , z_n)\in \mathfrak{H}^n$, 
$$ \gamma \cdot z \;\; = \;\; \bigg( \frac{a_1 z_1 + b_1}{c_1 z_1 + d_1}\; , \; \dots \;, \;  \frac{a_n z_n + b_n}{c_n z_n + d_n}\bigg) \;\;\;\;\;\;\;\; \text{ where } \;\; \sigma_j(\gamma) \; = \; \begin{pmatrix} a_j & b_j \\ c_j & d_j \end{pmatrix}.$$
We construct the Laplacian $\Delta$ on $SL_2(\mathfrak{o}) \backslash \mathfrak{H}^n$ from the usual Laplacians on the factors:
$$\Delta \;\; = \;\; \frac{\, 1 \, }{n} \; (\Delta_1 \; + \; \dots \; + \; \Delta_n) \;\;\;\;\; \text{ where } \;\;\;\;\; \Delta_j \; = \; y_j^2 \left( \frac{\partial^2}{\partial x_j^2} \; + \;  \frac{\partial^2}{\partial y_j^2}\right).$$
This is a nonpositive symmetric operator.  We parametrize the eigenvalues by $\lambda_w = w(w-1)$.  For $\lambda_w$ to be nonpositive real, we need $w \in \tfrac{1}{2} \, + \, i\R \; \cup \; [0,1]$.

Elements of global automorphic Sobolev spaces for $SL_2(\mathfrak{o})\backslash \mathfrak{H}^n$ have spectral expansions, in terms of an orthonormal basis $\{ F \}$ of spherical cusp forms, the constant automorphic form $1$, and the continuous family of Eisenstein series $E_{s,\chi}$  with $s$ on the critical line and $\chi$ an unramified grossencharacter:
$$E_{s, \chi}(z) \;\;= \;\; \sum_{\gamma \in P \cap SL_2(\mathfrak{o}) \backslash SL_2(\mathfrak{o})} \; \prod_{j = 1}^n \left(\mathrm{Im}(\sigma_j(\gamma) \cdot z_j)\right)^s \cdot \chi_j(\mathrm{Im}(\sigma_j(\gamma) \cdot z_j)),$$
where, as usual, $P$ is the standard parabolic subgroup of upper triangular matrices.  Fix a basepoint $z_o \in \mathfrak{H}^n$.  The automorphic delta distribution $\delta$ at $z_o$ has a spectral expansion:
$$\delta \;\; = \;\; \sum_F \; \overline{F}(z_o) \cdot F \;\; + \;\; \frac{1}{\langle 1, 1\rangle} \;\; + \;\; \sum_{\chi} \; \frac{1}{4\pi i} \int_{\frac{1}{2} + i \R} E_{1-s, \, \overline{\chi}}(z_o) \cdot E_{s, \chi} \; ds,$$
converging in a negatively indexed automorphic Sobolev space.  For $\mathrm{Re}(w) > \tfrac{1}{2}$, there is a unique solution $u_w$ to the automorphic differential equation  $(\Delta - \lambda_w) u_w \; = \; \delta$, and its spectral expansion, converging in an automorphic Sobolev space, is
$$u_w \;\; = \;\; \sum_F \; \frac{\overline{F}(z_o) \cdot F }{\lambda_F - \lambda_w} \;\; + \;\; \frac{1}{(\lambda_1 - \lambda_w)\langle 1, 1\rangle} \;\; + \;\; \sum_{\chi} \; \frac{1}{4\pi i} \int_{\frac{1}{2} + i \R} \frac{E_{1-s, \, \overline{\chi}}(z_o) \cdot E_{s, \chi}}{\lambda_{s, \chi} - \lambda_w} \; ds,$$
where $\lambda_F$, $\lambda_1$, and $\lambda_{s, \chi}$ denote the $\Delta$-eigenvalues of the Hilbert-Maass waveforms $F$, $1$, and $E_{s, \chi}$ occurring in the spectral expansion.  We emphasize that this is an equality of functions in an automorphic Sobolev space; it is not necessary to require pointwise equality.  Since the eigenvalues corresponding to the cuspidal and residual spectrum are discrete, elementary estimates ensure the meromorphic continuation, in $w$, of the cuspidal and residual part of the spectral expansion, as a Sobolev-space-valued function.  However, as we will show, the continuous part of the spectral expansion exhibits branching: for each nontrivial unramified grossencharacter $\chi$, the corresponding integral has \emph{two branch points} on the critical line.  

%As a function on the group $G = SL_2(k \otimes_\Q \R) \approx \left(SL_2(\R)\right)^n$, 
%$$E_{s, \chi}(g) \;\;= \;\; \sum_{\gamma \in P \cap SL_2(\mathfrak{o}) \backslash SL_2(\mathfrak{o})} \; \varphi_{\chi, s}(\gamma \cdot g), \;\;\;\; \varphi_{s,\chi}\left( \left(\begin{smallmatrix} a & \ast \\ 0 & d \end{smallmatrix} \right)\right) \;\; = \;\; |a/d|^s \cdot \chi(a/d),$$
%where, as usual, $P$ is the standard parabolic subgroup of upper triangular matrices.  

\begin{note*} In contrast, when $n=1$, i.e. $k = \Q$, there is no branching since the continuous part of the spectral expansion of the fundamental solution, 
$$u_w \;\; = \;\; \sum_F \; \frac{\overline{F}(z_o) \cdot F }{\lambda_F - \lambda_w} \;\; + \;\; \frac{1}{(\lambda_1 - \lambda_w)\langle 1, 1\rangle} \;\; + \;\;  \; \frac{1}{4\pi i} \int_{\frac{1}{2} + i \R} \frac{E_{1-s}(z_o) \cdot E_{s,}}{\lambda_{s} - \lambda_w} \; ds,$$
does not involve a sum over grossencharacters.  In Section \ref{min-para-Eis} we prove the analogous fact for $GL_3$.
  \end{note*}

Fix a grossencharacter $\chi$.  Take real parameters $t_\chi \, = \, (t_1, \dots, t_n)$ with $t_1 + \dots + t_n \; = \; 0$ such that 
$$\chi(\alpha) \;\; = \;\; \sigma_1(\alpha)^{i t_1} \dots \sigma_n(\alpha)^{i t_n},$$
where $\alpha \in (k \otimes_\Q \R)^\times$, and let 
$$\lVert t_\chi \rVert^2 \;\; = \;\; \tfrac{1}{n}\;(|t_1|^2 \, + \, \dots \, +\,  |t_n|^2).$$
We will show that the $\chi^{\text{th}}$ integral in the spectral expansion of $u_w$ admits a \emph{pathwise} meromorphic continuation to the complex plane with exactly two branch points: $\tfrac{1}{2} \, \pm \, i \, \lVert t_\chi \rVert$, when $\chi$ is nontrivial.

Let  $\mathcal{I}_\chi(w)$ denote the $\chi^{\text{th}}$ integral in the spectral expansion of $u_w$, as follows.
$$\mathcal{I}_\chi(w) \;\; = \;\; \int_{\frac{1}{2} + i \R} \; \frac{E_{1-s, \overline{\chi}}(z_o) \, E_{s, \chi}}{\lambda_{s, \chi} - \lambda_w} \; ds  \;\;\;\;\;\;\;\; (\mathrm{Re}(w) \, > \, \tfrac{1}{2})$$
This is a Sobolev-space-valued function of the complex parameter $w$, defined in a right half plane.

Writing the eigenvalue in terms of $s$ and $t_\chi$, 
$$\lambda_{s, \chi} \;\; = \;\; \frac{1}{n} \; \big( (s+ it_1)(s+ it_1 - 1) \, + \, \dots \, + \, (s+ it_n)(s+i t_n - 1) \big),$$
we can see that the integrand has poles when the following is satisfied.
\begin{eqnarray*}
\tfrac{1}{n} \; \big( (s+ it_1)(s+ it_1 - 1) \, + \, \dots \, + \, (s+ it_n)(s+i t_n - 1) \big) & = & w(w-1) \\
\tfrac{1}{n} \; \big( (s+ it_1)(s+ it_1 - 1) \, + \, \tfrac{1}{4} \; + \; \dots \; + \; (s+ it_n)(s+i t_n - 1) \, + \, \tfrac{1}{4} \big) & = &  w(w-1) \, + \, \tfrac{1}{4} \\
\tfrac{1}{n} \; \big((s-\tfrac{1}{2} \, + \, it_1)^2 \; + \; \dots \; + \; (s-\tfrac{1}{2} \, + \, it_n)^2\big) & = &  (w - \tfrac{1}{2})^2 
\end{eqnarray*}
Since $\sum t_j = 0$ and $\tfrac{1}{n} \sum t_j^2 = \lVert t_\chi \rVert^2$, we have
$$(s-\tfrac{1}{2})^2 \; - \;  \lVert t_\chi \rVert^2 \;\; = \;\;  (w - \tfrac{1}{2})^2 .$$
Thus the integrand has poles at
$$s \;\; = \;\; \tfrac{1}{2} \, \pm \, \sqrt{(w-\tfrac{1}{2})^2 \, + \, \lVert t_\chi \rVert^2}.$$

\begin{theorem} Let $\chi$ be a nontrivial grossencharacter and $\mathcal{I}_\chi(w)$ be the $\chi^\text{th}$ integral in the automorphic spectral decomposition of the fundamental solution $u_w$, as defined above.  Let $\gamma_1$ and $\gamma_2$ be $w$-paths in $\C$, each originating at a point $w_0$ in the right half plane $\mathrm{Re}(w) > \tfrac{1}{2}$, crossing the critical line once, and terminating at a point $w'_0$ in the left half plane $\mathrm{Re}(w) < \tfrac{1}{2}$, with $\gamma_1$ crossing the critical line at a height greater in magnitude than $\lVert t_\chi\rVert$ and $\gamma_2$ crossing at a height less in magnitude than $\lVert t_\chi\rVert$.  Then pathwise meromorphic continuations $\mathcal{I}_{\chi, \gamma_1}(w)$ and $\mathcal{I}_{\chi,\gamma_2}(w)$ of $\mathcal{I}_\chi(w)$ along the paths $\gamma_1$ and $\gamma_2$, respectively, differ by a term of moderate growth, namely by
$$\mathcal{I}_{\chi,\gamma_1}(w) \, - \, \mathcal{I}_{\chi, \gamma_2}(w) \;\; = \;\;\frac{4\pi i \; \cdot \; E_{1-s(\chi, w), \overline{\chi}}(z_o) \; \cdot \; E_{s(\chi, w), \chi}}{1 - 2 s(\chi, w)},$$
where $s(\chi, w)$ is defined as follows.  For fixed $w$ in $\mathrm{Re}(s)> \tfrac{1}{2}$,  $s(\chi, w)$ is the pole of the integrand of $\mathcal{I}_\chi(w)$ in $\mathrm{Re}(s)> \tfrac{1}{2}$.  As $w$ crosses the critical line,  $s(\chi, w)$ is defined by analytic continuation.

\begin{figure}
\centering
\includegraphics[scale=0.3]{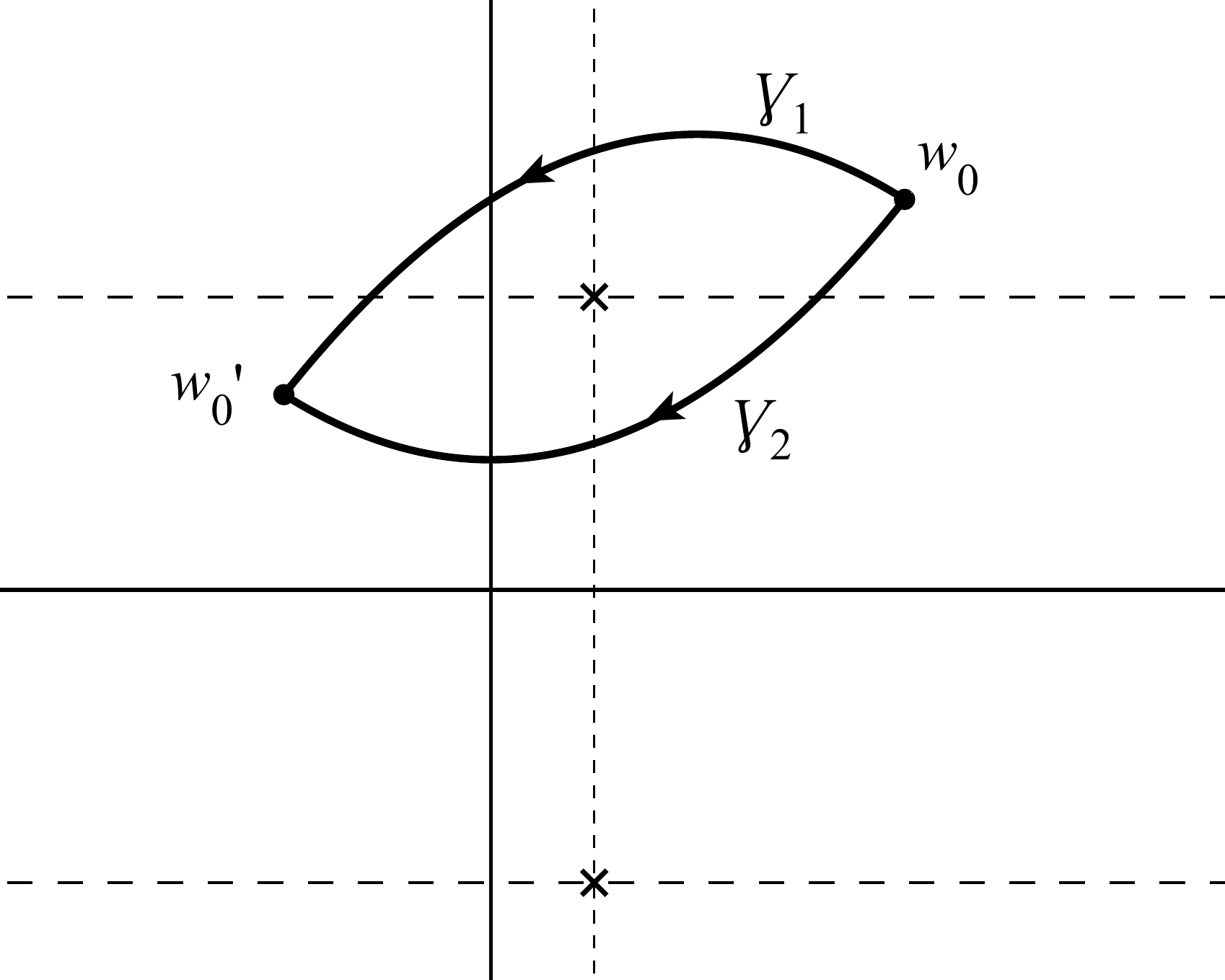}
\caption{Pathwise meromorphic continuation along these two paths in the $w$-plane yields functions that differ by a term of moderate growth.  The dotted vertical line is the critical line $\mathrm{Re}(w) = \tfrac{1}{2}$.  The dashed horizontal lines are $\mathrm{Im}(w) = \pm \lVert t_\chi\rVert$.}
\end{figure}

\begin{proof}
We meromorphically continue $\mathcal{I}_\chi(w)$ along two different paths.  For fixed $w$ to the right of the critical line, let $s(\chi, w)$ denote the pole of the integrand lying to the right of the critical line.  Since the numerator of the integrand is invariant under $s \to 1-s$, we regularize as follows.
\begin{eqnarray*}
\mathcal{I}_\chi(w) & = &   \int_{\frac{1}{2} + i \R} \; \frac{E_{1-s, \overline{\chi}}(z_o) \, E_{s, \chi} \; - \; E_{1-s(\chi,w), \overline{\chi}}(z_o) \, E_{s(\chi,w), \chi}}{\lambda_{s, \chi} - \lambda_w} \; ds\\
& & \;\;\;\;\;\;\;\; \;\; + \;\;  E_{1-s(\chi,w), \overline{\chi}}(z_o) \, E_{s(\chi,w), \chi} \; \cdot \; \int_{\frac{1}{2} + i \R}\; \frac{ds}{\lambda_{s, \chi} - \lambda_w} \;\;\;\;\;\;\;\; (\mathrm{Re}(w) \; > \; \tfrac{1}{2})
\end{eqnarray*}
By design the integrand of the first integral on the right side is continuous.  The second integral can be evaluated by residues.%, as long as the poles do not lie on the contour of integration, i.e. as long as $w$ does not lie on the critical line above or below the slit between $\tfrac{1}{2} \, \pm \, i \, \lVert t_\chi \rVert$.  Excluding this case, 
$$ \int_{\frac{1}{2} + i \R} \; \frac{ds}{\lambda_{s, \chi} - \lambda_w} \; = \; 2\pi i \, \times \, \operatorname*{Res}_{s \; = \;1 - s(\chi,w)} \; \frac{1}{(s - s(\chi,w))(s - (1-s(\chi,w))} \; = \;\frac{ 2\pi i}{1 - 2 s(\chi,w)} \;\;\;\;\;\;\;\; (\mathrm{Re}(w) \; > \; \tfrac{1}{2})$$
Thus we can rewrite the integral in the following way.
\begin{eqnarray*}
{\mathcal{I}_\chi(w)} & = &   \int_{\frac{1}{2} + i \R} \; \frac{E_{1-s, \overline{\chi}}(z_o) \, E_{s, \chi} \; - \; E_{1-s(\chi,w), \overline{\chi}}(z_o) \, E_{s(\chi,w), \chi}}{\lambda_{s, \chi} - \lambda_w} \; ds\\
& & \;\;\;\;\;\;\;\; \;\; + \;\;  E_{1-s(\chi,w), \overline{\chi}}(z_o) \, E_{s(\chi,w), \chi} \; \cdot \;\frac{ 2\pi i}{1 - 2 s(\chi,w)}  \;\;\;\;\;\;\;\; (\mathrm{Re}(w) \; > \; \tfrac{1}{2})
\end{eqnarray*}
We will see that when we move $w$ across the critical line, with imaginary part of greater magnitude than $\lVert t_\chi \rVert^2$, we pick up an Eisenstein series.  Note that since Eisenstein series do not lie in any finite-index automorphic Sobolev space, it is not trivial to determine where (in what function space) the (pathwise) meromorphic continuation lies.

Consider, for a moment, the case where $\chi=1$, the trivial character.  Then $s(\chi,w) \; = \; s(1,w) \; = \;  \tfrac{1}{2} \, + \, \sqrt{(w - \tfrac{1}{2})^2} \; = \; w$, since $\mathrm{Re}(w) \, > \,  \tfrac{1}{2}$ and $s(\chi, w)$ is defined to be to the right of the critical line.  Thus the integral corresponding to the trivial grossencharacter is as follows.
$$\mathcal{I}_1(w) \;\; = \;\;\int_{\frac{1}{2} + i \R} \; \frac{E_{1-s, 1}(z_o) \, E_{s, 1} \; - \; E_{1-w,1}(z_o) \, E_{w,1}}{\lambda_{s, 1} - \lambda_w} \; ds \;\; + \;\;  E_{1-w,1}(z_o) \, E_{w,1} \; \cdot \; \frac{ 2\pi i}{1 - 2 w} \;\;\;\;\;\;\;\; (\mathrm{Re}(w) \; > \; \tfrac{1}{2})$$
We move $w$ across the critical line and reverse the regularization.
$$\mathcal{I}_1(w) \; = \;\int_{\frac{1}{2} + i \R} \; \frac{E_{1-s, 1}(z_o) \, E_{s, 1}}{\lambda_{s, 1} - \lambda_w} \; ds \; - \; E_{1-w,1}(z_o) \, E_{w,1}\bigg( \int_{\frac{1}{2} + i \R} \; \frac{1}{\lambda_{s, 1} - \lambda_w} \; ds \; - \;  \frac{ 2\pi i}{1 - 2 w}\bigg) \;\;\;\;\;\;\; (\mathrm{Re}(w) \; < \; \tfrac{1}{2})$$
Now $s=w$ is the pole to the left of the critical line, and we evaluate the singular integral by residue calculus.
$$ \int_{\frac{1}{2} + i \R} \; \frac{ds}{\lambda_{s, \chi} - \lambda_w} \;\; = \;\; 2\pi i \, \times \, \operatorname*{Res}_{s \; = \; w} \;\; \frac{1}{(s - w)(s - (1-w))} \;\; = \;\;\frac{ 2\pi i}{2w-1} \;\;\;\;\;\;\;\; (\mathrm{Re}(w) \; < \; \tfrac{1}{2})$$
Thus we have the following expression for the integral corresponding to $\chi = 1$.
$$\mathcal{I}_1(w) \;\; = \;\;\int_{\frac{1}{2} + i \R} \; \frac{E_{1-s, 1}(z_o) \, E_{s, 1}}{\lambda_{s, 1} - \lambda_w} \; ds \;\; + \;\; E_{1-w,1}(z_o) \, E_{w,1} \; \cdot \;  \frac{ 4\pi i}{1-2w} \;\;\;\;\;\;\;\; (\mathrm{Re}(w) \; < \; \tfrac{1}{2})$$
From this, we see that the pathwise meromorphic continuation has an additional term, when $w$ is left of the critical line.

The pathwise meromorphic continuation also picks up an additional term in the case where $\chi$ is nontrivial, provided that $w$ crosses the critical line sufficiently far away from the real axis.  As $w$ crosses the critical line, with imaginary part greater in magnitude than $\lVert t_\chi \rVert$, the radicand, $(w - \tfrac{1}{2})^2 \; + \; \lVert t_\chi \rVert^2$, in the expression for $s(\chi,w)$ moves around the branch point of the square root, the origin.  Thus the analytic continuation of $s(\chi,w)$, along this path, is given as follows.
$$s(\chi,w) \;\; = \;\; \tfrac{1}{2} \, - \, \sqrt{(w-\tfrac{1}{2})^2 \, + \, \lVert t_\chi \rVert^2} \;\;\;\;\;\;\;\; (\mathrm{Re}(w) \, < \, \tfrac{1}{2})$$
As above, regularizing, evaluating the singular integral by residues, reversing the regularization, and again evaluating the singular integral, we obtain the meromorphic continuation of the $\chi^{\text{th}}$ integral, along a path $\gamma_1$ where $w$ crosses the critical line above a height of $\lVert t_\chi \rVert$ or below a height of $- \lVert t_\chi \rVert$.
$$\mathcal{I}_{\chi, \gamma_1}(w) \;\; = \;\;\int_{\frac{1}{2} + i \R} \; \frac{E_{1-s,\overline{\chi}}(z_o) \, E_{s, \chi}}{\lambda_{s, \chi} - \lambda_w} \; ds \;\; + \;\; E_{1-s(\chi, w),\overline{\chi}}(z_o) \, E_{s(\chi,w),\chi} \; \cdot \;  \frac{ 4\pi i}{1-2s(\chi,w)} \;\;\;\;\;\;\;\; (\mathrm{Re}(w) \; < \; \tfrac{1}{2})$$
Meromorphically continuing along a path $\gamma_2$ in which $w$ crosses the critical line with imaginary part within a distance of $\lVert t_\chi\rVert$ of the real axis does \emph{not} result in the additional Eisenstein series term, because, in this case, the radicand, $(w - \tfrac{1}{2})^2 \; + \; \lVert t_\chi \rVert^2$, in the expression for $s(\chi,w)$ stays strictly in the right half plane and thus does \emph{not} travel around the origin.  Thus branching is evident: pathwise meromorphic continuations of $\mathcal{I}_\chi(w)$ depend non-trivially on the path, the branch points being $w=\tfrac{1}{2} \pm i \, \lVert t_\chi \rVert$.  \end{proof}
\end{theorem}

\subsection{Additional Details: Branching of $s(\chi, w)$}

To understand better the different pathwise meromorphic continuations of $\mathcal{I}_\chi(w)$, the $\chi^{\text{th}}$ integral in the spectral expansion of the automorphic fundamental solution above, we explicitly parametrize $w$-paths crossing the critical line and show how the height of the crossing affects the radicand in the expression for the poles of the integrand.

For $\chi$ nontrivial, we can parametrize $w$ as follows
$$w \;\; = \;\; (\sigma + \tfrac{1}{2}) \; + \; (\alpha \, \lVert t_\chi \rVert) \, i,$$
with $\alpha \neq 0$.  To describe $w$ crossing the critical line on a horizontal path, we fix $\alpha$, and let $\sigma$ range from positive to negative values.  In terms of this parametrization, the radicand in the expression for the poles is
$$(w - \tfrac{1}{2})^2 \; + \; \lVert t_\chi\rVert^2 \;\; = \;\; \big(\sigma^2 \; + \; (1-\alpha^2) \lVert t_\chi \rVert^2\big) \; + \; \big(2\,\sigma\,\alpha\,\lVert t_\chi \rVert \big) \, i.$$
Let $x$ denote the real part of the radicand and $y$ the imaginary part.
\begin{eqnarray*}
 x & = & \sigma^2 \; + \; (1-\alpha^2) \lVert t_\chi \rVert^2\\ 
 y & = & 2\,\sigma \, \alpha \, \lVert t_\chi \rVert 
 \end{eqnarray*}
Eliminating the parameter $\sigma$, we can see that the curve is a right-facing parabola:
$$x \;\; = \;\; \frac{1}{4 \alpha^2 \, \lVert t_\chi \rVert^2} \; \big(y^2 \; + \; 4 \alpha^2(1-\alpha^2) \, \lVert t_\chi \rVert^4\big).$$
The direction that the radicand travels along this curve as $\sigma$ varies will depend on the sign of $\alpha$, i.e. it depends on whether $w$ is crossing below or above the real axis.  What is critical is whether the radicand travels around the origin, the branch point of the square root.  When $|\alpha| < 1$ (i.e. $w$ crosses the critical line between $\tfrac{1}{2} \pm i \, \lVert t_\chi \rVert$) the radicand does not travel around the origin, but when $|\alpha| > 1$, the radicand does travel around the origin.

%\begin{figure}[h!]
%\begin{minipage}[b]{0.5\linewidth}
%\centering
%\includegraphics[scale=0.25]{alpha0p6.pdf}
%\caption{$|\alpha|<1$}
%%\label{fig:figure1}
%\end{minipage}
%\hspace{0.5cm}
%\begin{minipage}[b]{0.5\linewidth}
%\centering
%\includegraphics[scale=0.25]{alpha1p4.pdf}
%\caption{$|\alpha| > 1$}
%%\label{fig:figure2}
%\end{minipage}
%\end{figure}
%

When $\chi$ is trivial, parametrize $w$ as $\; w \; = \; (\sigma + \tfrac{1}{2}) \; + \; i t_o$.  Fixing $t_o \neq 0$ and letting $\sigma$ range from positive to negative values describes $w$ crossing the critical line along a horizontal curve of height $t_o$.  In terms of these parameters, the radicand is
$$(w - \tfrac{1}{2})^2 \; + \; \lVert t_\chi\rVert^2 \;\; = \;\; \big(\sigma^2 \; - \; t_o^2\big) \; + \; \big(2\,\sigma\,t_o\big) \, i.$$
Again denoting the real and imaginary parts of the radicand by $x$ and $y$ respectively, we can see that the curve along which the radicand travels is
$$x \;\; = \;\; \frac{1}{4 t_o^2} \; (y - 2t_o^2)(y + 2 t_o^2).$$
This is a right-facing parabola, going around the origin.
%\begin{figure}[h!]
%\centering
%\includegraphics[scale=.25]{trivialchar.pdf}
%\caption{$\chi$ trivial}
%\end{figure}

The poles of the integrand of the $\chi^{\text{th}}$ term of the spectral expansion of $u_w$ are at
$$s \;\; = \;\; \tfrac{1}{2} \; \pm \; \sqrt{(w-\tfrac{1}{2})^2 \; + \; \lVert t_\chi \rVert^2}.$$
For fixed $w$ to the right of the critical line, we let $s(\chi, w)$ denote the pole to the right of the critical line.  We may choose a branch of the square root such that the following holds.
$$s(\chi, w) \;\; = \;\; \tfrac{1}{2} \; + \; \sqrt{(w-\tfrac{1}{2})^2 \; + \; \lVert t_\chi \rVert^2} \;\;\;\;\; (\mathrm{Re}(w) \; > \; \tfrac{1}{2})$$
As $w$ crosses the critical line, we analytically continue $s(\chi, w)$, and if the radicand travels around the origin, we have a sign change.
$$s(\chi, w) \;\; = \;\; \begin{cases} \tfrac{1}{2} \; + \; \sqrt{(w-\tfrac{1}{2})^2 \; + \; \lVert t_\chi \rVert^2} & (w\text{-path crosses sufficiently close to the real axis)}\\ \tfrac{1}{2} \; - \; \sqrt{(w-\tfrac{1}{2})^2 \; + \; \lVert t_\chi \rVert^2} & (w\text{-path crosses sufficiently far from the real axis)}\end{cases}$$

\section{Branching of $GL_3$ Automorphic Fundamental Solution}

Let $G = SL_3(\R)$, $K = SO(3)$ and $\Gamma = SL_3(\Z)$.  For simplicity we consider spherical automorphic forms.

Functions in global automorphic Sobolev spaces have spectral expansions, in terms of a spectral family of automorphic forms, consisting of cusp forms, Eisenstein series, and residues of Eisenstein series.  For $GL_3(\R)$, it suffices to take an orthonormal basis $\{ F \}$ of spherical cusp forms, the continuous family of minimal parabolic Eisenstein series $E^{1,1,1}_{\chi}$ where $\chi = \exp(\mu)$, for some $\mu \in \rho + i \mathfrak{a}^\ast$, and the family of $P^{2,1}$-Eisenstein series, $E^{2,1}_{f, s}$, with cuspidal data $f$ in an orthonormal basis of $GL_2$ cusp forms and complex parameter $s \in \tfrac{1}{2} + i \R$, along with the constant automorphic form (residue of minimal parabolic Eisenstein series).  In particular, for $\Phi$ in a $GL_3$ automorphic Sobolev space,
$$\Phi \;\; = \;\; \sum_{\text{cfm } F} \langle F, \Phi \rangle  \cdot F\; + \; \frac{\langle \Phi, 1 \rangle}{\langle 1,1 \rangle} \; + \; \frac{1}{|W|} \int_{\rho + i \mathfrak{a}^\ast} \langle \Phi , E^{1,1,1}_{\chi_\mu} \rangle \cdot E^{1,1,1}_{\chi_\mu} \, d \mu \; + \;
 \sum_{GL_2 \text{ cfms } f} \int_{\frac{1}{2} + i \R} \langle \Phi, E^{2,1}_{f,s} \rangle \cdot E^{2,1}_{f,s} \, ds,$$
where convergence is in a global Sobolev topology.  From now on, we drop the superscripts denoting the relevant parabolic for the Eisenstein series.

Fix a basepoint $x_0 \in G/K \approx \mathbb{H}^3$.  Then the automorphic delta distribution $\delta$ at $x_0$ has spectral expansion
$$\delta \;\; = \;\; \sum_{\text{cfm } F} \overline{F}(x_0) \cdot F \; + \; \frac{1}{\langle 1,1 \rangle} \; + \; \frac{1}{|W|} \int_{\rho + i \mathfrak{a}^\ast} E_{\bar{\chi}_\mu}(x_0) \cdot E_{\chi_\mu} \, d \mu \; + \;
\sum_{GL_2 \text{ cfms } f} \; \int_{\frac{1}{2} + i \R} E_{\bar{f},1-s}(x_0)  \cdot E_{f,s} \, ds,$$
converging in a negatively indexed automorphic Sobolev space.  For $\mathrm{Re}(w) > \tfrac{1}{2}$, there is a unique solution $u_w$ to the automorphic differential equation  $(\Delta - \lambda_w)^\nu \, u_w \; = \; \delta$, and its spectral expansion, converging in an automorphic Sobolev space, is
$$u_w \;\; = \;\; \sum_{\text{cfm } F}  \; \frac{\overline{F}(x_0)}{(\lambda_F - \lambda_w)^\nu} \, \cdot \, F \; + \; \frac{1}{\langle 1,1\rangle (\lambda_1 - \lambda_w)^\nu} \; + \; \frac{1}{|W|} \int_{\rho + i \mathfrak{a}^\ast} \; \frac{E_{\bar{\chi}_\mu}(x_0)}{(\lambda_{\chi} - \lambda_w)^\nu} \, \cdot \,  E_{\chi_\mu} \, d \mu $$
$$\hskip -3cm \; + \;  \sum_{GL_2 \text{ cfms } f} \;\;\; \int_{\frac{1}{2} + i \R} \; \frac{E_{\bar{f},1-s}(x_0)}{(\lambda_{f,s} - \lambda_w)^\nu} \cdot E_{f,s} \, ds,$$
where $\lambda_F$, $\lambda_1$, $\lambda_{\chi}$, and $\lambda_{f,s}$ denote the $\Delta$-eigenvalues of the waveforms $F$, $1$,  $E_{\chi}$, and $E_{f,s}$ occurring in the spectral expansion.  For simplicity, we choose $\nu$ to be the smallest integer, namely $\nu=2$,  that will ensure (using an automorphic Sobolev embedding theorem) uniform pointwise convergence of the spectral expansion.

Since the eigenvalues corresponding to the cuspidal and residual spectrum are discrete, elementary estimates ensure the meromorphic continuation, in $w$, of the cuspidal and residual part of the spectral expansion, as a Sobolev-space-valued function.  The part of the expansion corresponding to minimal parabolic Eisenstein series also admits a meromorphic continuation, as in the $GL_2$ case, but the part of the expansion corresponding to cuspidal data Eisenstein series exhibits branching: for each $GL_2$ cusp form $f$ in the chosen orthonormal basis, the corresponding integral has \emph{two branch points} on the critical line.  

\subsection{Pathwise meromorphic continuations of minimal parabolic Eisenstein series component} \label{min-para-Eis}

Let $\mathcal{I}_\chi(w)$ denote the most continuous part of the spectral expansion of the automorphic fundamental solution $u_w$.  
$$ \mathcal{I}_{\chi}(w)\;\; = \;\;  \int_{\rho + i \mathfrak{a}^\ast} \; \frac{E_{\bar{\chi}_\mu}(x_0)}{(\lambda_{\chi} - \lambda_w)^2} \, \cdot \,  E_{\chi_\mu} \, d \mu \;\;\;\;\;\;\;\; (\mathrm{Re}(w) > 0)$$
Then $\mathcal{I}_\chi(w)$ is a Sobolev-space-valued integral defined for $\mathrm{Re}(w) > 0$.  We show that $\mathcal{I}_\chi(w)$ does not exhibit branching in $w$.

Let $\mu = \rho + i \eta$, where $\eta \in \mathfrak{a}^\ast$.  Let $\lambda_w = w^2 - \Vert \rho \rVert^2$, where $w$ is a complex number with $\mathrm{Re}(w) > 0$.  Since the eigenvalue $\lambda_\chi$ is 
$$\langle \mu, \mu \rangle \, - \, 2\langle \mu, \rho \rangle \; = \; -(\langle \eta, \eta\rangle \, + \, \langle \rho, \rho \rangle) \;\; =\;\; -(\lVert \eta\rVert^2 \, + \, \lVert \rho \rVert^2),$$
the denominator of the integrand is $-(\rVert \eta\lVert^2 \, + \, w^2)$.  We rewrite the integral with these normalizations.
$$ \int_{\rho + i \mathfrak{a}^\ast} \; \frac{E_{\bar{\chi}_\mu}(x_0)}{(\lambda_{\chi} - \lambda_w)^2} \, \cdot \,  E_{\chi_\mu} \, d \mu \;\; = \;\; - \int_{ \mathfrak{a}^\ast} \; \frac{E_{\rho - i \eta}(x_0)}{(\lVert \eta \rVert^2 \, + \, w^2)^2} \, \cdot \,  E_{\rho + i \eta} \, d \eta$$
Thus the integrand is undefined when $w$ is purely imaginary and $\eta$ lies on the circle $\lVert \eta \rVert = | w|$ in $\mathfrak{a}^\ast \approx \R^2$.  Let $\mathcal{J}_w$ be the function-valued integral
$$\mathcal{J}_w \; = \; \int_{\lVert \eta \rVert \, = \, |w|} E_{\rho - i \eta}(x_0) \; E_{\rho + i \eta} \; d \eta.$$
Then we regularize.
$$\int_{ \mathfrak{a}^\ast} \; \frac{E_{\rho - i \eta}(x_0)}{(\lVert \eta \rVert \, + \, w^2)^2} \, \cdot \,  E_{\rho + i \eta} \, d \eta \;\; = \;\; \int_{ \mathfrak{a}^\ast} \; \frac{E_{\rho - i \eta}(x_0)\, \cdot \,  E_{\rho + i \eta}  \; - \; \mathcal{J}_w}{(\lVert \eta \rVert^2 \, + \, w^2)^2} \; d \eta \;\; + \;\; \mathcal{J}_w \cdot \int_{ \mathfrak{a}^\ast} \; \frac{1}{(\lVert \eta \rVert^2 \, + \, w^2)^2} \; d \eta$$
We evaluate the singular integral.
$$\int_{ \mathfrak{a}^\ast} \; \frac{1}{(\lVert \eta \rVert^2 \, + \, w^2)^2} \; d \eta \;\; = \;\; \int_{\R^2}\; \frac{1}{(\lVert \eta \rVert^2 \, + \, w^2)^2} \; d \eta \;\; = \;\;2\pi \; \cdot \; \int_0^\infty \frac{r \, dr}{(r^2 + w^2)^2} \;\; = \;\; \frac{\pi}{w^2}$$
So we have the following.
$$\int_{ \mathfrak{a}^\ast} \; \frac{E_{\rho - i \eta}(x_0)}{(\lVert \eta \rVert \, + \, w^2)^2} \, \cdot \,  E_{\rho + i \eta} \, d \eta \;\; = \;\; \int_{ \mathfrak{a}^\ast} \; \frac{E_{\rho - i \eta}(x_0)\, \cdot \,  E_{\rho + i \eta}  \; - \; \mathcal{J}_w}{(\lVert \eta \rVert^2 \, + \, w^2)^2} \; d \eta \;\; + \;\; \frac{\pi \; \mathcal{J}_w }{w^2}  \;\;\;\;\; (\mathrm{Re}(w) > 0)$$

Now we may move $w$ across the imaginary axis and undo the regularization.  Since the value of the singular integral is again $\pi/w^2$, the extra terms cancel, as folllows.

\begin{eqnarray*}
\lefteqn{\int_{ \mathfrak{a}^\ast} \; \frac{E_{\rho - i \eta}(x_0) \, \cdot \,  E_{\rho + i \eta} }{(\lVert \eta \rVert \, + \, w^2)^2} \; d \eta}  \\
& & \\
& = &  \int_{ \mathfrak{a}^\ast} \; \frac{E_{\rho - i \eta}(x_0)\, \cdot \,  E_{\rho + i \eta}}{(\lVert \eta \rVert^2 \, + \, w^2)^2} \; d \eta   \;\; - \;\; \mathcal{J}_w \cdot \int_{ \mathfrak{a}^\ast} \; \frac{1}{(\lVert \eta \rVert^2 \, + \, w^2)^2} \; d \eta\;\; + \;\; \frac{\pi \; \mathcal{J}_w }{w^2}  \;\;\;\;\; (\mathrm{Re}(w) <0)\\
& & \\
& = &  \int_{ \mathfrak{a}^\ast} \; \frac{E_{\rho - i \eta}(x_0)\, \cdot \,  E_{\rho + i \eta}}{(\lVert \eta \rVert^2 \, + \, w^2)^2} \; d \eta   \;\; - \;\; \frac{\pi \; \mathcal{J}_w }{w^2}  \;\; + \;\; \frac{\pi \; \mathcal{J}_w }{w^2}  \;\;\;\;\; (\mathrm{Re}(w) <0)\\
%& & \\
%& = &  \int_{ \mathfrak{a}^\ast} \; \frac{E_{\rho - i \eta}(x_0)\, \cdot \,  E_{\rho + i \eta}}{(\lVert \eta \rVert^2 \, + \, w^2)^2} \; d \eta    \;\;\;\;\; (\mathrm{Re}(w) <0)
\end{eqnarray*}

\subsection{Pathwise meromorphic continuations of cuspidal Eisenstein series component}

Let $f$ be a $GL_2$ cusp form occuring in the orthonormal basis chosen above, and let $s_f(s_f-1)$ be its eigenvalue, with $s_f \in [0,1] \cup \tfrac{1}{2} + i \R$.  If we let $s_f = \tfrac{1}{2} + i t_f$, then $t_f \in -i[0, \tfrac{1}{2}] \cup \R$, with $t_f = -i/2$ corresponding to $s_f = 0, 1$.  Let $\mathcal{I}_f(w)$ denote the integral corresponding to $f$ in the spectral expansion of $u_w$.
$$\mathcal{I}_f(w)\;\; = \;\; \int_{\frac{1}{2} + i \R} \; \frac{E_{\bar{f},1-s}(x_o) \, E_{f, s}}{(\lambda_{f,s} - \lambda_w)^2} \;\; ds  \;\;\;\;\;\;\;\; (\mathrm{Re}(w) \, > \, \tfrac{1}{2})$$
This is a Sobolev-space-valued integral defined in a right half plane.  

The cuspidal data Eisenstein series $E_{f,s}$ generates a principal series, induced from character on the minimal parabolic $P=MN$, which is determined by its action on $M$:
$$\left(\begin{smallmatrix} a_1 & & \\ & a_2 & \\ & & a_3 \end{smallmatrix}\right) \;\; \mapsto \;\; |a_1|^{s_f + s} \; |a_2|^{-s_f + s} \; |a_3|^{-2s}.$$
As derived in the appendix, the eigenvalue of Casimir on a minimal parabolic Eisenstein series $E_\chi$ is
$$\lambda_{\chi} \;\; = \;\; 2( s_1^2  \, + \, s_1 s_2 \, + \, s_2^2 \, - \, 2 s_1 \, - \, s_2), \;\;\;\;\;\;\;\; \text{where } \;\;  \chi \left(\begin{smallmatrix} a_1 & & \\ & a_2 & \\ & & a_3 \end{smallmatrix}\right) \;\; = \;\; |a_1|^{s_1} \; |a_2|^{s_2} \; |a_3|^{s_3}.$$
Thus the eigenvalue of Casimir on $E_{f,s}$ is
$$\lambda_{f,s} \;\; = \;\; 2\big( s_f(s_f-1) \, + \, 3s(s-1)\big).$$
Letting $\lambda_w = 6w(w-1)$, we can see that the integrand has poles when the following is satisfied.
\begin{eqnarray*}
2\big( s_f(s_f-1) \, + \, 3s(s-1)\big) & = & 6 w(w-1) \\
\big((s_f - \tfrac{1}{2})^2 - \tfrac{1}{4}\big) \; + \; 3\big((s - \tfrac{1}{2})^2 - \tfrac{1}{4}\big) & = &  3\big((w- \tfrac{1}{2})^2  \, - \, \tfrac{1}{4}\big) \\
(s_f - \tfrac{1}{2})^2 - \tfrac{1}{4}\; + \; 3(s - \tfrac{1}{2})^2 & = & 3 (w- \tfrac{1}{2})^2\\
(s - \tfrac{1}{2})^2 & = &  (w- \tfrac{1}{2})^2 \;-\; \tfrac{1}{3}\big((s_f - \tfrac{1}{2})^2 - \tfrac{1}{4}\big)\\
(s - \tfrac{1}{2})^2 & = &  (w- \tfrac{1}{2})^2 \;+\; \tfrac{1}{3}\big(t_f^2 + \tfrac{1}{4}\big)\\
\end{eqnarray*}
Thus the integrand has poles at
$$s \;\; = \;\; \tfrac{1}{2} \, \pm \,  \, \sqrt{(w-\tfrac{1}{2})^2 \, + \tfrac{1}{3}(t_f^2 + \tfrac{1}{4})}.$$

\begin{theorem} Let $f$ be a spherical cusp form in the chosen basis of $GL_2$ cusp forms and $\mathcal{I}_f(w)$ be the corresponding integral in the automorphic spectral decomposition of the fundamental solution $u_w$, as defined above.  Let $\gamma_1$ and $\gamma_2$ be $w$-paths in $\C$, each originating at a point $w_0$ in the right half plane $\mathrm{Re}(w) > \tfrac{1}{2}$, crossing the critical line once, and terminating at a point $w'_0$ in the left half plane $\mathrm{Re}(w) < \tfrac{1}{2}$, with $\gamma_1$ crossing the critical line at a height greater in magnitude than $\sqrt{\tfrac{1}{3}(t_f^2 + \tfrac{1}{4})}$ and $\gamma_2$ crossing at a height less in magnitude than $\sqrt{\tfrac{1}{3}(t_f^2 + \tfrac{1}{4})}$.  Then pathwise meromorphic continuations $\mathcal{I}_{f,\gamma_1}(w)$ and $\mathcal{I}_{f,\gamma_2}(w)$ of $\mathcal{I}_f(w)$ along the paths $\gamma_1$ and $\gamma_2$ respectively differ by a term of moderate growth, namely by
$$\mathcal{I}_{f, \gamma_1}(w) \, - \, \mathcal{I}_{f, \gamma_2}(w) \;\; = \;\; \frac{8\pi i \; \cdot \; E_{\bar{f}, 1-s(f, w)}(x_o) \; \cdot \; E_{f, s(f, w)}}{(1 - 2 s(f, w))^3},$$
where $s(f, w)$ is defined as follows.  For fixed $w$ in $\mathrm{Re}(s)> \tfrac{1}{2}$,  $s(f, w)$ is the pole of the integrand of $\mathcal{I}_f(w)$ in $\mathrm{Re}(s)> \tfrac{1}{2}$.  As $w$ crosses the critical line,  $s(f, w)$ is defined by analytic continuation.

\begin{proof} We meromorphically continue $\mathcal{I}_f(w)$ along two different paths.  For fixed $w$ to the right of the critical line, let $s(f, w)$ denote the pole of the integrand lying to the right of the critical line.  Since the numerator of the integrand is invariant under $s \to 1-s$, we regularize as follows.
\begin{eqnarray*}
{\mathcal{I}_f(w)} & = &   \int_{\frac{1}{2} + i \R} \; \frac{E_{\bar{f}, 1-s}(x_0) \, E_{f, s} \; - \; E_{\bar{f}, 1-s(f,w)}(x_0) \, E_{f, s(f,w)}}{(\lambda_{f,s} - \lambda_w)^2} \; ds\\
& & \;\;\;\;\;\;\;\; \;\; + \;\;  E_{\bar{f}, 1-s(f,w)}(x_0) \, E_{f, s(f,w)} \; \cdot \; \int_{\frac{1}{2} + i \R}\; \frac{ds}{(\lambda_{f,s} - \lambda_w)^2} \;\;\;\;\;\;\;\; (\mathrm{Re}(w) \; > \; \tfrac{1}{2})
\end{eqnarray*}
By design the integrand of the first integral on the right side is continuous.  The second integral can be evaluated by residues.%, as long as the poles do not lie on the contour of integration, i.e. as long as $w$ does not lie on the critical line above or below the slit between $\tfrac{1}{2} \, \pm \, i \, \lVert t_\chi \rVert$.  Excluding this case, 
$$ \int_{\frac{1}{2} + i \R} \; \frac{ds}{(\lambda_{f, s} - \lambda_w)^2}  =  2\pi i \operatorname*{Res}_{s \; = \;1 - s(f,w)} \; \frac{1}{(s - s(f,w))^2(s - (1-s(f,w))^2}  = \frac{ 4\pi i}{(2 s(f,w)-1)^3} \;\;\;\;\;\; (\mathrm{Re}(w) \; > \; \tfrac{1}{2})$$
Thus we can rewrite $\mathcal{I}_f(w)$.
\begin{eqnarray*}
{\mathcal{I}_f(w)} & = &   \int_{\frac{1}{2} + i \R} \; \frac{E_{\bar{f}, 1-s}(x_0) \, E_{f, s} \; - \; E_{\bar{f}, 1-s(f,w)}(x_0) \, E_{f, s(f,w)}}{(\lambda_{f,s} - \lambda_w)^2} \; ds\\
& & \;\;\;\;\;\;\;\; \;\; + \;\;  E_{\bar{f}, 1-s(f,w)}(x_0) \, E_{f, s(f,w)} \; \cdot \; \frac{ 4\pi i}{(2 s(f,w)-1)^3} \;\;\;\;\;\;\;\; (\mathrm{Re}(w) \; > \; \tfrac{1}{2})
\end{eqnarray*}

As $w$ crosses the critical line, with imaginary part greater in magnitude than $\sqrt{\tfrac{1}{3}(t_f^2 + \tfrac{1}{4})}$, the radicand, $(w - \tfrac{1}{2})^2 \; + \; \tfrac{1}{3}(t_f^2 + \tfrac{1}{4})$, in the expression for $s(f, w)$ moves around the branch point of the square root, the origin.  Thus the analytic continuation of $s(f, w)$, along this path, is as follows.
$$s(f, w) \;\; = \;\; \tfrac{1}{2} \, - \, \sqrt{(w-\tfrac{1}{2})^2 \, + \, \tfrac{1}{3}(t_f^2 + \tfrac{1}{4})} \;\;\;\;\;\;\;\; (\mathrm{Re}(w) \, < \, \tfrac{1}{2})$$
Regularizing, evaluating the singular integral by residues, reversing the regularization, and again evaluating the singular integral, we obtain the meromorphic continuation of $\mathcal{I}_f(w)$, along a path $\gamma_1$ where $w$ crosses the critical line above a height of $\sqrt{\tfrac{1}{3}(t_f^2 + \tfrac{1}{4})}$ or below a height of $- \sqrt{\tfrac{1}{3}(t_f^2 + \tfrac{1}{4})}$.
$$\mathcal{I}_{f, \gamma_1}(w) \;\;= \;\; \int_{\frac{1}{2} + i \R} \; \frac{E_{\bar{f}, 1-s}(x_0) \, E_{f, s}}{\lambda_{f, s} - \lambda_w} \; ds \;\; + \;\; E_{\bar{f}, 1-s(f, w)}(x_0) \, E_{f, s(f, w)} \; \cdot \;  \frac{ 8\pi i}{(1-2s(f, w))^3} \;\;\;\;\;\;\;\; (\mathrm{Re}(w) \; < \; \tfrac{1}{2})$$
Meromorphically continuing along a path $\gamma_2$ in which $w$ crosses the critical line with imaginary part within a distance of $\sqrt{\tfrac{1}{3}(t_f^2 + \tfrac{1}{4})}$ of the real axis does \emph{not} result in the additional Eisenstein series term, because, in this case, the radicand, $(w - \tfrac{1}{2})^2 \; + \; \tfrac{1}{3}(t_f^2 + \tfrac{1}{4})$, in the expression for $s(f, w)$ stays strictly in the right half plane and thus does \emph{not} travel around the origin.  Thus branching is evident: pathwise meromorphic continuations of $\mathcal{I}_f(w)$ depend non-trivially on the path, the branch points being $\tfrac{1}{2} \pm i \, \sqrt{\tfrac{1}{3}(t_f^2 + \tfrac{1}{4})}$.\end{proof}
\end{theorem}

\section{Appendix: Eigenvalue of Casimir on minimal parabolic Eisenstein series}

The eigenvalue of Casimir on a minimal parabolic Eisenstein series with character $e^\mu$, $\mu = \rho + i \eta$, where $\rho$ is half the sum of positive roots and $\eta \in \mathfrak{a}^\ast$, is
$$\langle \mu, \mu \rangle - 2 \langle \mu, \rho \rangle \;\; = \;\; - \big( \langle \eta, \eta \rangle \, + \, \langle \rho, \rho \rangle\big).$$
Consider $\mu \in \mathfrak{a}_\C^\ast$ as a linear combination of positive roots with coefficients $s_\alpha$, $s_\beta$, and $s_{\alpha + \beta}$:
$$\mu \;\; = \;\; s_\alpha \, \alpha \; + \; s_\beta \, \beta \; + \; s_{\alpha + \beta} \, (\alpha + \beta) \;\; = \;\; (s_\alpha + s_{\alpha + \beta}) \, \alpha \; + \; (s_\beta + s_{\alpha + \beta}) \, \beta.$$
Then
$$e^\mu\left(\begin{smallmatrix} a_1 & & \\ & a_2 & \\ & & a_3 \end{smallmatrix}\right) \;\; = \;\; \left|\frac{a_1}{a_2}\right|^{s_\alpha}  \cdot \;  \left|\frac{a_2}{a_3}\right|^{s_\beta}  \cdot \;  \left|\frac{a_1}{a_3}\right|^{s_{\alpha+\beta}}\;\; = \;\;  \left|\frac{a_1}{a_2}\right|^{s_\alpha + s_{\alpha + \beta}} \cdot \; \left|\frac{a_2}{a_3}\right|^{s_\beta + s_{\alpha + \beta}}$$
and the eigenvalue is
$$2 \big( s_{\alpha}^2 \, + \, s_\beta^2 \, -\,  s_\alpha s_\beta \, + \, s_\alpha s_{\alpha + \beta} \, + \, s_\beta s_{\alpha + \beta} \, -\, s_\alpha \, - \, s_\beta \, - \, 2s_{\alpha + \beta}\big).$$
On the other hand, it is also common to parametrize the character as
$$e^\mu\left(\begin{smallmatrix} a_1 & & \\ & a_2 & \\ & & a_3 \end{smallmatrix}\right) \;\; = \;\; |a_1|^{s_1}  \cdot \; |a_2|^{s_2} \cdot \; |a_3|^{s_3},\;\;\;\;\;\;\;\; \text{ where } \;\; s_1 + s_2 + s_3 \; = \; 0.$$
In this case, $\mu \; = \; s_1 \alpha \, + \, (s_1 + s_2) \beta$, so the eigenvalue is
$$2 \big( s_1^2 \, + \, s_1 s_2 \, + \, s_2^2 \, - \, 2s_1 \, - \, s_2\big).$$

%% BIBLIOGRAPHY %%%%%%%%%

\end{document}